\documentclass[final,3p]{elsarticle}
\usepackage{amssymb}
 \usepackage{amsthm}
\usepackage{amsmath,bbm,dsfont, color}

\journal{arXiv}

\begin{document}

\newcommand*{\id}{\mathds{1}}
\newcommand*{\zahlen}{\mathbb{Z}}
\newcommand*{\en}{\mathbb{N}}
\newcommand*{\er}{\mathbb{R}}
\newcommand{\esssup}{\mathop{\text{ess\:sup}}}
\newcommand{\htimes}{\mathop{\text{\large$�$}}}
\newcommand{\hexists}{\mathop{\text{\LARGE$\exists$}}}
\newcommand{\hforall}{\mathop{\text{\LARGE$\forall$}}}
\newtheorem{lem}{Lemma}
\newtheorem{theo}{Theorem}
\newtheorem{cor}{Corollary}
\newtheorem{prop}{Proposition}
\newtheorem{mydef}{Definition}
\newtheorem{rem}{Remark}

\begin{frontmatter}

 \title{Does the fully parabolic quasilinear 1D Keller-Segel system enjoy long-time asymptotics analogous to its parabolic-elliptic simplification?}

\author{Jan Burczak \fnref{fn1}} \ead{jb@impan.pl} \address{Institute of Mathematics, Polish Academy of Sciences, \'Sniadeckich 8, 00-950 Warsaw.}

\author{Tomasz Cie\'slak \fnref{fn2}} \ead{T.Cieslak@impan.pl} \address{
Institut f\"{u}r Mathematik, Universit\"{a}t zu Z\"{u}rich, Winterthurerstrasse 190, 8057 Z\"{u}rich,\\
Institute of Mathematics, Polish Academy of Sciences, \'Sniadeckich 8, 00-950 Warsaw.
}

\author{Cristian Morales-Rodrigo \fnref{fn3}} \ead{cristianm@us.es} \address{Dpto. de Ecuaciones
Diferenciales y An\'alisis Num\'erico, Fac. de Matem\'aticas, Univ. de Sevilla,
Calle Tarfia s/n, 41012-Sevilla, Spain}

\fntext[fn1]{supported by the International Ph.D. Projects Programme of Foundation for Polish Science within the Innovative Economy Operational Programme 2007-2013}

 \fntext[fn2]{partially supported by the Polish Ministry of Science and Higher Education under grant
number NN 201366937} \fntext[fn3]{supported by Ministerio de Ciencia e Innovaci\'on and FEDER under grant MTM2009-12367}


\begin{abstract}
We show that the one-dimensional fully parabolic Keller-Segel system with nonlinear diffusion possesses global-in-time
solutions, provided the nonlinear diffusion is equal to $ \frac{1}{(1+u)^\alpha}, \; \alpha<1, $ independently on the
volume of the initial data. We also show that in the critical case, i.e. for $\alpha = 1$, the same result holds for
initial masses smaller than a prescribed constant. Additionally, we prove existence of initial data for which solution
blows up in a finite time for any nonlinear diffusion integrable at infinity. Thus we generalize the known blowup result
of parabolic-elliptic case to the fully parabolic one. However, in the parabolic-elliptic case the above mentioned
integrability condition on nonlinear diffusion sharply distinguishes between global existence and blowup cases.
We are unable to recover the entire global existence counterpart of this result in a fully parabolic case. 
\end{abstract}

\begin{keyword}
 fully parabolic Keller-Segel system \sep global existence \sep  finite-time blowup

\end{keyword}




\end{frontmatter}

\section{Introduction}
\noindent
The Keller-Segel model was introduced as a system of four parabolic quasilinear equations in \cite{KS:1} to describe in a
mathematical way the motion of small organisms under the chemotactic forces. One of its issues was to provide a model
whose solution aggregates in a finite-time as a result of activity of chemotactic force attracting cells. Moreover this
aggregation was supposed to be caused by large enough initial mass. This gives an argument to claim that spontaneous
self-organization is possible in nature and it could be simply a consequence of laws of physics. This could be seen as a
particular example of differentiation of cells leading to a better evolutionary adaptation.  In \cite{Nan}
Nanjundiah simplified the original model introducing the so-called minimal version of the Keller-Segel model consisting
of two parabolic equations. Basing on numerical experiments it was claimed that the simplified model is still a proper
description of an aggregation phenomenon. Moreover the attention was also paid to the threshold value of initial mass
guaranteeing the finite-time blowup of solution, a phenomenon interpreted as aggregation of cells. However,
rigorous results that followed were showing finite-time blowups of solutions merely of further parabolic-elliptic
simplifications of a model (see: \cite{jl:expl}, \cite{Na95}, \cite{Nagai2}). Nagai and collaborators ( \cite{Na95},
\cite{Nagai2}, \cite{NSY97}) have  even found the values of initial mass yielding finite-time blowups  both for radial
solutions in a ball and nonradial ones in more general two-dimensional domains. In the radial case the conjecture in
\cite{child} based on numerical computations  was confirmed. Let us emphasize that the models were parabolic-elliptic
ones, though there is a wide agreement in the community that model with two parabolic equations is the one describing
the reality better. In \cite{HV:1} the authors pointed one solution blowing up in a fully parabolic case. Unluckily,
no exact information on a threshold value which distinguishes between finite-time blowup and infinite lifespan of
solution is available. There is no mathematical proof that results from parabolic-elliptic case hold also in a fully
parabolic one. Recently, in \cite{Racz}, there was presented an argument matching behaviour of solutions in
parabolic-elliptic case with fully parabolic one. However, this result shows the stability properties of a model when
passing from a fully parabolic to the parabolic-elliptic case in a very wide class of functions. Due to this obstacle
one cannot say if the finite-time blowup from a parabolic-elliptic case is inherited in a fully parabolic one.

On the other hand, recently in \cite{clKS} it was proved that solutions to the one-dimensional quasilinear fully
parabolic case blow up provided nonlinear diffusion is weak enough. Next in \cite{clSP} and \cite{clGlob_vs_blow} 
the exact strength of the diffusion distinguishing between finite-time blowup and global solutions in a corresponding parabolic-elliptic one-dimensional systems was identified. Our aim in the present paper is to study if the same results are available in a fully parabolic case. The confirmation would be a nice evidence that the results in a parabolic-elliptic case can be believed to hold also in a fully parabolic one.

In a present paper we are unfortunately unable to realize this program to a full extent. This is a first step of studying
the connection between one-dimensional parabolic-elliptic and fully parabolic models. Actually, with respect to the class
of nonlinear diffusions admitting finite-time blowup in a fully parabolic case we prove an analogous result as in
\cite{clGlob_vs_blow}, where the parabolic-elliptic system was studied. However our global-in-time existence result is
weaker than its parabolic-elliptic counterpart. Namely, we prove that global solutions exist for the subcritical
diffusions (which is already a new result in the fully parabolic case and fulfills the global existence counterpart to
the finite-time blowup result in \cite{clKS}) without any restriction on the size of the initial mass, but the result is
not pushed far enough to cover all the nonlinearities nonintegrable at infinity. Basing on our result one cannot say what
happens if one starts with an arbitrarily large initial mass in the case of the nonlinearity which could be a candidate
for a critical one, namely $a(u)=\frac{1}{1+u}$. In this case we prove a global existence for initial mass smaller than
a certain threshold. The question about behaviour of solutions emanating from masses larger than our threshold remains
an open problem. So is the question of qualitative behaviour of solutions to the Keller-Segel system ran by
diffusions which are nonintegrable at infinity, but weaker than $\frac{1}{1+u}$. 

Let us emphasize that studying the lifespan of solutions in a fully parabolic case requires completely different methods than those used in a parabolic-elliptic one. A change of variable introduced in \cite{clSP} which reduces a parabollic-elliptic system to one equation posessing a Liapunov functional is not available in a fully parabolic case. That reformulation essentially simplifies further studies of properties of solutions to the parabolic-elliptic Keller-Segel system in one dimension.

Let us now describe our results in a more precise way.

We consider a following one dimensional Keller-Segel problem with nonlinear difussion $a$
\begin{multline*}[KS]
\begin{cases}
u_t = (a (u) u_x - {\chi} u v_x)_x \quad in \; (0, \infty) \times (0,1),\; \chi>0, \\
\varepsilon v_t = v_{xx}  + u - M \quad in \; (0, T) \times (0,1), \; \varepsilon>0, \\
u_x =  v_x = 0  \quad on \; (0, T) \times \{0;1\}  \\
(u, v) (0) = (u_0, v_0) .
\end{cases}
\label{1_ks}
\end{multline*}
We save the letters $u, v $ to denote the solution of the above [KS] system and in the entirety of this paper we
usually refer to solution of [KS] system simply as to $u, v$. We denote absolute value of a real number as
$| \cdot |$, $L_p $ norm as $| \cdot |_p$ and $W^{1,p}$ norm as $| \cdot |_{1,p}$. $C, K$ are constants, which may
change even in a single line. \vskip 2mm
\noindent
This paper is divided into four main sections and organized as follows: In the remainder of this introduction we collect
well-known facts on short-time existence, uniqueness and regularity of solutions as well as availability of the Liapunov
functional. The next two sections deal with dependence of time-asymptotics of solutions to [KS] on  nonlinear
diffusion. The second one is devoted to the problem of global-in-time existence of $L^\infty$-bounded solutions of [KS].
In the third part, we present finite-time blowup result for supercritical (integrable) diffusions $a$. Finally, we conclude with some remarks on possible extensions of our results and its connections to some functional inequalities.

\noindent
Let us state the standard well-posedness and classical solvability result
\begin{prop}[Well-posedness, classical solvability, conservation of mass for [KS{]} system]
Assume that
\begin{equation}
a \in C^2(\er_+), \; a> 0, \quad M := \int_0^1 u_0 (t,x) dx > 0, \; u_0 \ge 0,  \quad \int_0^1 v_0 (t,x) dx = 0,  \quad u_0, v_0 \in W^{1,2} (0,1)
\end{equation}
 in [KS]. Then it admits local-in-time, unique, classical solution with maximal time of existence $T_m$
 \begin{equation}
 (u, v) \in \; C ( [0, T_m) \times [0,1]; \er^2) \; \cap \; C^{2,1} ( (0, T_m) \times [0,1]; \er^2),
 \end{equation}
  which additionally satisfies
  \begin{equation}\label{pp}
  \int_0^1 u (t,x) dx= \int_0^1 u_0  (t,x) dx = M, \; u> 0\; \; \mbox{for}\; \; t>0,  \quad \int_0^1 v (t,x) dx = 0.
 \end{equation}
 Moreover, if $T_m < \infty$ we have blowup of $L^\infty $ norm, i.e.
 \begin{equation}
 \lim_{t \rightarrow T_m} [ |u(t)|_\infty + |v(t)|_\infty ] \rightarrow \infty.
 \end{equation}
\label{prop1.1}
\end{prop}
\noindent
The proof follows the general theory of parabolic (triangular) systems. For precise references consult \cite[p.440]{clKS}.\\
Interestingly, the [KS] system posesses the Liapunov functional, which, in one-dimensional setting, is additionally bounded from below. Let us introduce:
\begin{mydef} For $b: \er_+ \rightarrow \er$ such that: $ b'' (s) = \frac{a(s)}{s}, \;  b' (1) = 0, \; b(1) = 0$ and $b(u_0) \in L^1$, where $a$ denotes diffusion associated with [KS],  let
\begin{equation}\label{Lapun}
\lambda (t) : = \int_0^1 \left( b(u(t,x)) - u(t,x)v(t,x) + \frac{1}{2} |v_x (t,x)|^2 dx \right)
\end{equation}
 for $(u,v)$ solving [KS].
 \label{def1.1}
\end{mydef}
\noindent
\begin{prop}
Under assumptions of Proposition \ref{prop1.1} it holds for $ t \in [0, T_m) $:
\begin{equation}\label{Logr}
\frac{d}{dt} \lambda (t)  \le -\varepsilon |v_t|^2_2 (t), \qquad \lambda (t)  \ge -\frac{M^2}{2}.
\end{equation}
\label{prop1.2}
\end{prop}
\noindent
This proposition is shown in \cite[Lemmas 4, 5]{clKS}.

Next we have a Corollary which is a simple consequence of (\ref{Lapun}) and (\ref{Logr}).
\begin{cor}\label{cor}
Under assumptions of Proposition \ref{prop1.1} there exists $C>0$ such that
\begin{equation}\label{ogr}
\int_0^{T_m}\int_0^1 |v_t|^2 dxdt< C.
\end{equation}
\end{cor}

\section{Global existence}
\noindent
The main result of this sections reads
\begin{theo}
For $(u, v)$ solving [KS] it holds:
\begin{equation}
|u(t)|_\infty \le C
\end{equation}
 provided one of the following assumptions is valid:
\begin{itemize}
\item[i.] (subcritical case) nonlinear diffusion takes the form: $a(s) = \frac{1}{(1+s)^\alpha }$ for $\alpha \in [0,1)$ and the initial mass $\int_0^1 u_0(x) dx = M$ is arbitrary;
\item[ii.] (critical case) nonlinear diffusion takes the form: $a(s) = \frac{1}{(1+s) }$ and the initial mass $\int_0^1 u_0(x)dx = M$ satisfies $M< \frac{2}{\sqrt{\chi}} - 1$.
\end{itemize}
\label{glo}
\end{theo}
\noindent
For the sake of clarity, first we derive the main ingredients needed to show the above theorem and in the final part of
this section we combine them in the main proof. Let us begin by quoting the following regularity result:
\begin{prop}[boundedness of highly integrable nonnegative solutions to quasilinear parabolic equation]
Consider the following system
\begin{multline}
\begin{cases}
u_t =  {\rm div} (D (x, t, u) \nabla u + f (x,t)) \quad in \; (0, T) \times (0,1) \\
 \partial_\nu v = 0  \quad on \; (0, T) \times \{0,1\}
\end{cases}
\label{wt}
\end{multline}
Assume that:
\begin{enumerate}
\item[(1)] $D \in C^1 ([0,1] \times [0, T) \times [0, \infty) )$, $D>0$ and there exist $m \in \er, s_0 \ge 0, \delta > 0$ such that $D(x, t, s) \ge \delta s^{m-1} $ for $s \ge s_0$;
\item[(2)] $f \in C^0 (0,T; C^0 (\bar{\Omega}) \cap C^1 (\Omega) )$, $f_x(0)=f_x(1)=0$ for  $t\in (0, T)$ , $f \in L^\infty (L^{q_1}) $;
\item[(3)] $ u \ge 0$ solves \eqref{wt} and $u \in L^\infty (L^{p_0}) $;
\item[(4)] following inequalities for $q_1, p_0$ hold:
\begin{equation}
q_1 > 3
\label{wt4.1}
\end{equation}
\begin{equation}
p_0 \ge 1, \quad p_0 > 1 - m \frac{2 q_1 - 3}{q_1 - 3},  \quad  p_0 > 1-m,   \quad  p_0 > \frac{1-m}{2}
\label{wt4.2}
\end{equation}
\end{enumerate}
then
\begin{equation}
|u(t)|_\infty \le C.
\end{equation}
\label{glo.5}
\end{prop}
\begin{proof}
This theorem is a version of  \cite[Lemma A.1]{winkler-tao}. More precisely, taking $n=1, g \equiv 0$  in \cite[Lemma A.1]{winkler-tao} one obtains exactly Proposition \ref{glo.5}.
\end{proof}
\noindent
Later on we need the following maximal regularity result for the one-dimensional heat equation in Sobolev spaces. The formulation which we
need involves a right hand-side in $L^\infty(0,T;L^1(0,1))$. For the proof we refer to \cite[Lemma 4.1]{horstmann-winkler}.
\begin{prop}[maximal regularity for one-dimensional heat equation]
\label{prop1.3}
Let $v\in C^{1,2}([0,T)\times [0,1])$ be a solution to the heat equation
\begin{equation}\label{rownanie}
v_t-v_{xx}=f \;\mbox{in}\;(0,T)\times (0,1) \;\;,v(0,x)=v_0(x)\in C[0,1],
\end{equation}
where $f\in L^\infty(0,T;L^1(0,1))\cap C([0,T)\times[0,1])$. Moreover $v$ satisfies boundary conditions
\begin{equation}\label{warbrz}
v_x(0)=v_x(1)=0.
\end{equation}
Then
\[
\sup_{t\in[0,T)}|v_x|_{q}<C
\]
for any $q<\infty$.
\end{prop}

The following lemma is the main ingredient of the proof of Theorem \ref{glo} in the subcritical diffusion case.
\begin{lem}\label{lemm}
Assume that nonlinear diffusion of [KS] $a(u) = \frac{1}{(1+u)^\alpha}$ with $\alpha \in [0,1]$. Additionally, assume
that solution $ (u, v) $ of [KS] satisfies for time independent constant $C$
\begin{equation}
|(1+u)^{\frac{s- \alpha +1 }{2}}|_p (t) \le C,
\label{2.1.A1}
\end{equation}
where $s, p$ are positive numbers satisfying
\begin{equation}
p (s- \alpha +1 ) > 2 \alpha, \quad p \ge 1,
\label{2.1.A2}
\end{equation}
then
\begin{equation}
|(1+u)|_{s+1} (t) \le C(u_0, v_0, p, \alpha) .
\label{2.1.the}
\end{equation}
\label{prop2.1}
\end{lem}
\begin{proof}
Without loss of generality set $\varepsilon =\chi= 1$ in [KS]. Multiplying the first equation of [KS] by $(1+u)^s$,
integrating over space and performing one integration by parts we have
\begin{equation*}
\frac{1}{s+1} \frac{d}{dt} \int  (u+1)^{s+1} = - s \int a(u) |(u+1)_x|^2 (u+1)^{s-1} + s \int (u+1-1)v_x (u+1)_x   (u+1)^{s-1},
\end{equation*}
which by assumption $a(u) = \frac{1}{(1+u)^\alpha}$ gives
\begin{multline*}
\frac{d}{dt} \int  (u+1)^{s+1}  + (s+1) s \int |(u+1)_x|^2 (u+1)^{s-1-\alpha} = \\
 s (s+1) \int v_x (u+1)_x   (u+1)^{s} - s (s+1) \int v_x (u+1)_x   (u+1)^{s-1} = \\
s (s+1) \int [ (u+1)_x   (u+1)^{\frac{s-1-\alpha}{2}}]   (u+1)^{\frac{s+1+\alpha}{2}} v_x
- s (s+1) \int [ (u+1)_x   (u+1)^{\frac{s-1-\alpha}{2}}]   (u+1)^{\frac{s-1+\alpha}{2}} v_x.
\end{multline*}
Hence,
\begin{multline*}
\frac{d}{dt} \int  (u+1)^{s+1}  + \frac{4s(s+1)}{2 (s-\alpha + 1)^2} \int |[(u+1)^{\frac{s-\alpha + 1}{2}}]_x|^2  \le
\frac{s(s+1)}{2}  \int |v_x|^2 [ (u+1)^{s+1+\alpha} + (1+u)^{s-1+\alpha} ],
\end{multline*}
adding to both sides of the above inequality the term $\int (u+1)^{s+1}$ we arrive at
\begin{multline}
\frac{d}{dt} \int  (u+1)^{s+1} + \int (u+1)^{s+1} + \frac{2 s (s+1)}{(s-\alpha + 1)^2} \int |[(u+1)^{\frac{s-\alpha + 1}{2}}]_x|^2  \le \\
\frac{s(s+1)}{2}  \int |v_x|^2 [ (u+1)^{s+1+\alpha} + (1+u)^{s-1+\alpha} ]  + \int (u+1)^{s+1}.
\label{2.1.1}
\end{multline}
To proceed further, recall that the Gagliardo-Nirenberg interpolation inequality in one-dimensional domains yields
\begin{equation*}
|w|_r^r \le K |w|^{r \theta}_{1,2} |w|_p^{r(1- \theta)}, \quad \; \frac{1}{r} = - \frac{\theta}{2} + \frac{1- \theta}{p}.
\end{equation*}
Therefore, requiring $r \theta = 2$ one obtains $r=2p+2$. Setting $ w : = (u+1)^{\frac{s- \alpha +1}{2}}$ we have
\begin{equation}
\int (u+1)^{(2p+2) \frac{s- \alpha +1}{2}} \le K |(u+1)^{\frac{s- \alpha +1}{2}}|^{2}_{1,2} |(u+1)^{\frac{s- \alpha +1}{2}}|_p^{2p}
\label{2.1.2}
\end{equation}
In order to use this inequality, we add  $ \frac{2 s (s+1)}{(s-\alpha + 1)^2} \int (u+1)^{s - \alpha +1}$ to both sides
of \eqref{2.1.1}, thus
\begin{multline*}
\frac{d}{dt} \int  (u+1)^{s+1} + \int (u+1)^{s+1} + \frac{2 s (s+1)}{(s-\alpha + 1)^2}  |(u+1)^{\frac{s-\alpha + 1}{2}}|_{1,2}^2  \le \\
 \frac{2 s (s+1)}{(s-\alpha + 1)^2} \int (u+1)^{s - \alpha +1} + \frac{s(s+1)}{2}  \int |v_x|^2 [ (u+1)^{s+1+\alpha} + (1+u)^{s-1+\alpha} ]  + \int (u+1)^{s+1} \\
 \le \left(1 + \frac{2 (s+1)}{s} \right)\int (u+1)^{s+1} +  \frac{2 (s+1)}{s}  + \frac{s(s+1)}{2} \int  |v_x|^2 [ 1 + 2(u+1)^{s+1+\alpha} ]
\end{multline*}
where the second inequality holds by $ cx \le c + cx^\gamma, \; \gamma \ge 1, \: c, x $ positive and $\alpha \le 1$.
In view of assumption \eqref{2.1.A1} the above inequality and Young's inequality imply for any positive, fixed smallness
constants $\eta, \delta$
\begin{multline*}
\frac{d}{dt} \int  (u+1)^{s+1} + \int (u+1)^{s+1} +
\left(K |(u+1)^{\frac{s- \alpha +1}{2}}|^{2}_{1,2} |(u+1)^{\frac{s- \alpha +1}{2}}|_p^{2p}\right)
\frac{2 s (s+1)}{K (s-\alpha + 1)^2}  |(u+1)^{\frac{s- \alpha +1}{2}}|_p^{-2p}\le \\
 \frac{3s+2}{s} \int (u+1)^{s+1}  + \delta  \int (u+1)^{s+1+\alpha + \eta}   + \frac{\eta}{s+1 + \alpha + \eta} \left(\delta \frac{s+1 + \alpha + \eta}{s + 1 + \alpha} \right)^{ -\frac{s+1 + \alpha}{\eta}}
\int  |v_x|^{\frac{ 2(s+1+\alpha + \eta)}{\eta}} \\+ \frac{s(s+1)}{2}\int  |v_x|^2  +  \frac{2 (s+1)}{s} .
\end{multline*}
Next, by \eqref{2.1.A1} and \eqref{2.1.2}
\[
\frac{d}{dt} \int  (u+1)^{s+1} + \int (u+1)^{s+1} +
\left(\int (u+1)^{(p+1) (s- \alpha +1)}\right)
\frac{2 s (s+1)}{K (s-\alpha + 1)^2}  C^{-2p}\le
\]
\[
2 \delta  \int (u+1)^{s+1+\alpha + \eta}   + \frac{\eta}{s+1 + \alpha + \eta} \left(\delta \frac{s+1 + \alpha + \eta}{s + 1 + \alpha} \right)^{ -\frac{s+1 + \alpha}{\eta}} \int  |v_x|^{\frac{2( s+1+\alpha + \eta)}{\eta}} +
\]
\begin{equation}\label {2.1.4}
\frac{s(s+1)}{2}\int  |v_x|^2  +  \frac{2 (s+1)}{s} + \frac{3s+2}{s} \left( \frac{(3s+2))s+1)}{\delta s (s+1+ \alpha + \eta)}\right)^\frac{s+1}{\alpha + \eta}.
\end{equation}
Taking $\eta, \delta$ such that
\[
2 \delta = \frac{2 s (s+1)}{K (s-\alpha + 1)^2}  C^{-2p}
\]
and
\[
(p+1) (s- \alpha +1) = s+1+\alpha + \eta,
\]
the latter is possible due to assumption  \eqref{2.1.A2}  which is equivalent to
\[
(p+1) (s- \alpha +1) > s+1+\alpha,
\]
and abandoning a precise control over constants, we finally arrive  at
\begin{equation}
\frac{d}{dt} \int  (u+1)^{s+1} + \int (u+1)^{s+1} \le C(s, \alpha) \left(1+ |v_x|_q^q \right)
\end{equation}
for some number $q = q(s,\alpha)$.

Thus in view of Proposition \ref{prop1.3} we have global-in-time bounds for $\int  (u+1)^{s+1}$.
\end{proof}
\noindent
At this point we will apply an iterative procedure, using at each step the previous proposition. The method is
implemented in the following corollaries. First we consider the subcritical case which allows us to obtain global-in-time boundedness of arbitrarily large $L_p$ norm of $u$.
\begin{cor}
Let $u, v$ solve the [KS] system with subcritical diffusion $a(u) = (1+u)^{-\alpha}, \; \alpha \in [0,1)$, then
\begin{equation}
|u|_p (t) \le C(u_0, v_0, p, \alpha)
\end{equation}
for any $p \in [1, \infty)$, where $C(u_0, v_0, p, \alpha)$ is time independent.
\label{cor2.2}
\end{cor}
\begin{proof}
We choose $s_i, p_i; \: i=1,2,..$ as follows
\begin{equation}
s_1 = 1+ \alpha, \qquad
p_1 = \frac{2}{s_1 - \alpha +1 } = 1
\end{equation}
\begin{equation}
s_{i+1} =  2s_i + 1+ \alpha , \qquad
p_{i+1}  = \frac{2 (s_i + 1)}{s_{i+1} - \alpha +1 } = 1
\end{equation}
which implies
\begin{equation}
s_{i} =  (1+ \alpha) (2^i -1), \qquad p_i (s_i - \alpha +1 ) \ge 2 > 2 \alpha.
\end{equation}
Therefore such choice of $s_i, p_i$ satisfies assumption \eqref{2.1.A2} for every $i$. Now we recursively obtain that
\begin{enumerate}
\item for $s_1, p_1$ from Proposition \ref{prop1.1} we have
\[
|(u+1)^{\frac{s_1-\alpha +1}{2}}|_{p_1} = \int (u+1) = M +1
\]
 so assumption \eqref{2.1.A1} holds and by Lemma \ref{lemm} we obtain $|u|_{s_1 +1 } (t) \le C$;
\item  for $s_{i+1}, p_{i+1}$ we get:
\[
|(u+1)^{\frac{s_{i+1}-\alpha +1}{2}}|_{p_{i+1}}  = \int (u+1)^{s_i +1},
\]
so assumption \eqref{2.1.A1} is valid by virtue of inductive assumption and consequently via Lemma \ref{lemm} one obtains:
\[
|u|_{s_{i+1} +1 } (t) \le C.
\]
\end{enumerate}
 Hence we lift  stepwise the integrability of $u$ to any fixed number $p<\infty$.
\end{proof}
\noindent
Observe that for the case of critical diffusion the bound $|u|_1 (t) \leq C$ is not enough to apply Lemma 1.
The reason is that \eqref{2.1.A1} holds for $p_1 = 2/s_1$, which in turn violates \eqref{2.1.A2}.
Nevertheless, we can
obtain the following weaker result for the critical diffusion.
\begin{cor}
Let $u, v$ solve the [KS] system with critical diffusion $a(u) = (1+u)^{-1}$, then for any $p \in [1, \infty)$
\begin{equation}
|u|_{1+ \varepsilon} (t) \le C \;\;\mbox{implies}\;\; |u|_p (t) \le C(u_0, v_0, p, \alpha)
\end{equation}
where $C(u_0, v_0, p, \alpha)$ is time independent and $\varepsilon> 0$ is arbitrarily small.
\label{cor2.3}
\end{cor}
\begin{proof}
Choose $s_i, p_i, \: i=1,2,..\;$ as follows
\begin{equation}
s_1 = 2(1+\varepsilon),  \qquad
p_1 = \frac{2(1+\varepsilon) }{s_1 } = 1,
\end{equation}
\begin{equation}
s_{i+1} =  2s_i + 2,  \qquad
p_{i+1}  = \frac{2 (s_i + 1)}{s_{i+1} } = 1.
\end{equation}
Again we can proceed inductively using Lemma \ref{lemm} to show that
\[
|u|_p (t) \le C(u_0, v_0, p, \alpha)
\]
for any $p<\infty$. We see that
\begin{itemize}
\item[(i=1)] $ |(u+1)^\frac{s_1}{2}|_{p_1} = | (u+1)^{1+ \varepsilon}|_1 \le C, \quad p_1 s_1 = 2+ 2 \varepsilon >2 = 2 \alpha$ by definitions of $s_1, p_1$ and assumption $|u|_{1+ \varepsilon} (t) \le C $;
\item[(i+1)] $ |(u+1)^\frac{s_{i+1}}{2}|_{p_{i+1}} = | (u+1)^{1+ s_i}|_1 \le C, \quad p_{i+1} s_{i+1} = 2+ 2 s_i >2 = 2 \alpha$ by definitions of $s_i, p_i$ and recursive assumption.
\end{itemize}
\end{proof}
\noindent
At this stage, we possess all the needed ingredients for showing Theorem \ref{glo} in the subcritical case, because Corollary \ref{cor2.2} allows us to obtain high integrability and Proposition \ref{glo.5} enables us to perform the step from high integrability to boundedness. However,  in the critical diffusion case, we lack the bound on $ |u|_{1+ \varepsilon} (t) $. In what follows we struggle to obtain one. Firstly we derive global-in-time bounds for $|u|_{LlogL}$ under assumptions on smallness of M (see Lemma \ref{lem2.5}). Next we utilize an idea from \cite{debye} to perform the step from the bound $|u|_{LlogL} (t) \le C$ to $|u|_{1+ \varepsilon} (t) \le C$.\\
In order to show Lemma  \ref{lem2.5} we need the following result.
\begin{prop}
For a function $m\in W^{1,2} (0,1)$, it holds for every $\nu > 0$
\begin{equation}
\int e^{2m} \leq \frac{1 + \nu}{4} \left(\int e^m \right)^2 \int |m_x|^2 +\left(1 + \frac{1}{ \nu} \right)  \left(\int e^m \right)^2.
\label{2.4.1}
\end{equation}
\label{prop2.4}
\end{prop}
\begin{proof}
First, for arbitrary $m \in W^{1,1} (0,1)$ a constant in the one-dimensional Sobolev imbedding is $1$, as for any
$x \in (0,1)$ holds
\begin{equation}
|m(x)| \le |m (x) - m(z)| +|m(z)| \le \int_0^1 |m_x| + |m(z)| \le |m|_{1,1}.
\end{equation}
Hence
\[
\int e^{2m} \leq \left(\sup_{x \in (0,1)} e^{\frac{m(x)}{2}} \right)^2 \int e^m \le   |e^{\frac{m}{2}}|_{1,1}^{2}  \int e^m =\left( \frac{1}{2} \int | e^{\frac{m}{2}} m_x| +  \int e^{\frac{m}{2}} \right)^2  \int e^m .
\]
Next, by H\"{o}lder's inequality and $ (a+b)^2 \le (1+ \nu) a^2 + (1 + \frac{1}{ \nu} ) b^2 $ we have
\begin{multline*}
\int e^{2m}\leq
\left[ \frac{1}{2} \left(\int e^m\right)^{\frac{1}{2}}\left(\int |m_x|^2\right)^{\frac{1}{2}} +  \int e^{\frac{m}{2}} \right]^2  \int e^m
\le
\left(\int e^m\right)^2 \left[  \frac{1 + \nu}{4} \left(\int |m_x|^2\right) + \left(1 + \frac{1}{ \nu} \right) \right].
\end{multline*}
\end{proof}
\begin{cor}\label{corlog}
Let $u$ be a solution to [KS] with $ \int_0^1 u_0(x)dx = M$. For $ M  < \frac{2}{\sqrt{\chi}}-1$ the
following inequality holds
\begin{equation}\label{imp}
\chi \int_0^1 u^2 dx < \int_0^1\left|(\log (1+u))_x\right|^2dx + C(M, \chi).
\end{equation}
\end{cor}
\begin{proof}
Due to (\ref{pp}), $u$ is a positive function with fixed mass. Letting  $m:=\log(1+u)$, i.e. $u+1=e^m$, in the
Proposition \ref{prop2.4} leads to
\[
\int_0^1u^2dx < (M+1)^2 \frac{1 + \nu}{4} \int_0^1\left|(\log (1+u))_x\right|^2dx + C(M, \nu),
\]
which gives thesis provided
\[
\frac{1 + \nu}{4}\chi (M+1)^2  < 1
\]
holds. In view of Proposition \ref{prop2.4} one can take $\nu $
as small as needed, thus the assumed inequality $ M  < \frac{2}{\sqrt{\chi}}-1$ is equivalent to
$\chi (M+1)^2 \frac{1}{4} < 1$.
\end{proof}
Next lemma is inspired by the additional Liapunov functional for Keller-Segel system that appeared in \cite{GZ:ch}.
\begin{lem}
For $u, v$ solving the [KS] system with critical diffusion $a(u) = (1+u)^{-1}$, such that
$ M = \int_0^1 u_0(x)dx <\frac{2}{\sqrt{\chi}}-1$, there is a global bound for $LlogL $ norm of
$u+1$, i.e. $\int |(u+1)log(u+1)| (t) \le C$.
\label{lem2.5}
\end{lem}
\begin{proof}
We test first equation of [KS] system with $log(u+1)$ and the second one twice, once with
$\chi(u+1)$ and another time with $- \chi log(u+1)$, obtaining
\begin{equation*}
\begin{array}{cccc}
\displaystyle\int (u+1)_t log (u+1) &=& -\displaystyle \int \frac{|(u+1)_x|^2}{(u+1)^2} + \chi\int \frac{u}{u+1}v_x (u+1)_x \\[3mm]
&=& -\displaystyle\int \frac{|(u+1)_x|^2}{(u+1)^2} + \chi\int v_x (u+1)_x-\chi\int\frac{1}{u+1}v_x (u+1)_x,\\[3mm]
\varepsilon \chi \displaystyle \int v_t (u+1)& =& - \chi \displaystyle \int  v_x (u+1)_x + \chi \int u (u+1) - \chi M \int (u+1),\\[3mm]
- \varepsilon \chi \displaystyle \int v_t log (u+1)& =& +\chi\displaystyle\int\frac{1}{u+1}v_x (u+1)_x - \chi \int u log(u+1) + \chi M \int log (u+1).
\end{array}
\end{equation*}
Adding the above equalities we arrive at
\[
\frac{d}{dt} \int (u+1) log (u+1) + \int |(log(u+1))_x|^2 =
\]
\[
=\chi\int u (u+1)-\chi\int u\log(u+1) -\chi M (M+1)+\chi M\int \log(u+1) - \chi\varepsilon \int v_t ((u+1)-\log(u+1)),
\]
Since for any $\mu>0$ we can find such $C$ that
\[
\int (u+1) log(u+1) = \int  |(u+1) log(u+1)| \le \mu\int u^2 + C\;\;\mbox{and}\;\; \log(u+1)\leq u\;\mbox{for}\;u \ge 0,
\]
Young's inequality yields
\begin{multline}
\frac{d}{dt} \int (u+1) log (u+1) + \int (u+1) log (u+1) +  \int |(log(u+1))_x|^2 \le \\
\chi\int u^2 +\chi\varepsilon \int |v_t| (2u+1)+ \mu\int u^2 + C\le (\chi+2\mu)\int u^2 + C|v_t|_2^2 + C.
\label{2.5.4}
\end{multline}
Thanks to \eqref{imp} and arbitrary smallness of $\mu$ \eqref{2.5.4} yields
\[
\frac{d}{dt} \int (u+1) log (u+1) + \int (u+1) log (u+1) \le C|v_t|_2^2 + C(M, \chi),
\]
which in turn gives boundedness of $ \left[ \int_0^1 (u(t,x)+1) log (u(t,x)+1)dx \right]$ in view of Corollary \ref{cor}.
\end{proof}
\begin{rem}\label{uwaga}
Observe that we need a smallness assumption on  $M$ in Lemma \ref{lem2.5} only due to restrictions of 
Proposition \ref{prop2.4} and Corollary \ref{corlog}. One sees that in order to reach Lemma \ref{lem2.5} without any restriction 
on the size of initial mass or $\chi$ it would suffice to show the following claim:
for any $\delta > 0$ and any function $m\in W^{1,2}(0,1)$ such that 
\[
\int_0^1 e^m dx=M,
\]
the following inequality holds
\begin{equation}\label{sehrwichtig}
\int e^{2m} < \delta M^2\int |m_x|^2 + h_\delta \left( \int e^m \right), 
\end{equation}
where $h_\delta$ denotes an arbitrary family of continuous functions $h_\delta: [0,\infty)\rightarrow \er$. 
Unluckily the claim is false; see Lemma \ref{nieprawda} in the last section of this paper.
\end{rem}
\noindent

Now we state an inequality which we need in the sequel. It is actually a one-dimensional version of the inequality in \cite{debye}. The proof is exactly the same, we present it only for reader's convenience.
\begin{prop}
For $w \in W^{1,2} (0,1)$ and any $\delta>0$ there exists $C_{\delta}$ such that the following inequality holds
\begin{equation}
|w|_4^4 \le \delta |w|_{1,2}^2 |w|_{LlogL} + C_{\delta} |w|_1.
\label{2.7.0}
\end{equation}
\label{prop2.7}
\end{prop}
\begin{proof}
Define $\eta_N : \er \rightarrow \er$ as follows:
\begin{equation}
\eta_N (s) = \left\{ \begin{array}{clcr}
0 & for \; |s| \le N, \\
2(|s| - N) & for \: |s| \in (N, 2N],\\
|s| & for \; |s| > 2N.
\end{array} \right.
\label{2.7.1}
\end{equation}
The Gagliardo-Nirenberg inequality gives
\begin{equation}
|\eta (w)|_4^4 \le K | \eta(w)|_{1,2}^2 |\eta(w)|_1^2.
\label{2.7.2}
\end{equation}
We estimate right-hand side of \eqref{2.7.2} in the
 following manner
\begin{equation}
\begin{array}{clcl}
|\eta(w)|_{1,2}^2 &=& \displaystyle\int |\eta' (w) w_x|^2 + \int |\eta (w)|^2 \le 4 |w|_{1,2}^2, \\
|\eta(w)|_1 &\le& \displaystyle\int |w| \id_{ \{ |w| > N \} }  = \int |wlogw| (logw)^{-1} \id_{ \{ |w| > N \} } \le
 (logN)^{-1} |w|_{LlogL},
\end{array}
\label{2.7.3}
\end{equation}
where in the first equation we used $|\eta'| \le 2$ and in the second we take $N > e$. Next we estimate the difference between $w$ and $\eta(w)$
\begin{equation}
||w|-\eta(w)|_4^4 = \int ||w|-\eta(w)|^4  \id_{ \{ |w| \le 2N \} } \le (2N)^3 \int ||w|-\eta(w)|  \id_{ \{ |w| \le 2 N \} } \le 8N^3 |w|_1^1.
\label{2.7.5}
\end{equation}
Considering \eqref{2.7.3}, \eqref{2.7.5} we can write in view of \eqref{2.7.2}
\[
|w|_4^4 \le 8 [||w|-\eta(w)|_4^4 + |\eta (w)|_4^4]  \le  64 N^3 |w|_1^1 + 8K | \eta(w)|_{1,2}^2 |\eta(w)|_1^2 \le 64 N^3 |w|_1^1 + 32 K  |w|_{1,2}^2  (logN)^{-2} |w|^2_{LlogL},
\]
which for $N = \max \left\{ e; e^{\sqrt{\frac{32K}{\delta}}} \right\}$ is \eqref{2.7.0}.
\end{proof}
\noindent
The next lemma follows the method introduced in \cite{NSY97}.
\begin{lem}
If $u$ solving [KS] with critical diffusion $a(u) = (1+u)^{-1}$ admits $|u+1|_{LlogL} (t) \le C$, then $|u+1|_3 (t) \le C$.
\label{lem2.8}
\end{lem}
\begin{proof}
Without loss of generality we assume $\chi=1$. Next we test the first equation in [KS] with $(u+1)^2$ to get
\[
\frac{1}{3} \frac{d}{dt} \int (u+1)^3 = -2 \int |(u+1)_x|^2 + 2\int (u+1)^2 u_x v_x-2\int (u+1) u_xv_x
\]
\[
=-2 \int |(u+1)_x|^2 - 2/3\int (u+1)^3v_{xx}+\int (u+1)^2v_{xx}
\]
\[
=-2 \int |(u+1)_x|^2-2/3\int (u+1)^3(v_t+M-u)+\int (u+1)^2(v_t+M-u)
\]
\[
=-2 \int |(u+1)_x|^2  -\frac{2M}{3}\int (u+1)^3 -2/3\int(1+u)^3v_t  + \int(1+u)^2v_t  +2/3\int u (u+1)^3 - \int u (1+u)^2 + M \int (u+1)^2.
\]
Thus for $w:=u+1$
\begin{equation}\label{wyjsciowe}
\frac{1}{3} \frac{d}{dt} \int w^3+ \frac{2M}{3} \int w^3+2 \int |w_x|^2\leq 5/3\int w^3 |v_t|+ \int w^4+ \frac{3 M^2}{4}.
\end{equation}
Recall that in view of (\ref{2.7.0}) for arbitrary small $\delta>0$ holds
\begin{equation}\label{oszac1}
\int w^4 \le \delta \left[ \int |w_x|^2 + \int w^2 \right] \int |w logw| + C_{\delta} \int w.
\end{equation}
On the other hand by the Gagliardo-Nirenberg inequality used twice
\[
 |w|_6^3 \leq C |w|_{1,2}^{\frac{3}{5}} |w|_3^{\frac{12}{5}} =  C |w|_{1,2}^{\frac{3}{5}} |w|_3^{\frac{3}{2}} |w|_3^{\frac{9}{10}} \le C |w|_{1, 2}^{\frac{3}{5}} |w|_3^{\frac{3}{2}}|w|_{1,2}^{\frac{2}{5}}|w|_1^\frac{1}{2} = C(M) |w|_{1, 2} |w|_3^{\frac{3}{2}},
\]
which means for any $\theta>0$
\begin{equation}\label{oszac2}
 \int w^3|v_t| \leq |v_t|_2 |w|_6^3 \leq \theta|w_x|_2^2+ \theta|w|_2^2 + C_{\theta}|v_t|_2^2|w|_3^3.
\end{equation}
In view of (\ref{oszac1}) and (\ref{oszac2}) and due to assumption $|w|_{LlogL} (t) \le C$ , (\ref{wyjsciowe}) implies that $f(t):=\int_0^1w^3(x,t)dx$ satisfies
\begin{equation}\label{glowne}
\frac{d}{dt}f(t)+Mf(t)\leq C |v_t|_2^2 f(t)+C,
\end{equation}
which equivalently can be rewritten 
\[
\frac{d}{dt} f(t) +\left(M-C  |v_t|_2^2\right) f(t)\leq C.
\]
Hence
\begin{equation}\label{NG}
f(t) \le f(0)e^{-\phi(t)}+C\left(\int_0^t e^{\phi(s)} ds\right)e^{-\phi(t)}, 
\end{equation}
where
\[
\phi(t)=\int_0^t \left(M-C|v_t(s)|_2^2\right) ds.
\]
Taking into account \eqref{ogr} we have
\[
-\phi(t)=-Mt+C\int_0^t|v_t|_2^2\leq -Mt+C
\]
and 
\[
\left(\int_0^t e^{\phi(s)} ds\right)e^{-\phi(t)}\leq \frac{e^{Mt}}{M}e^{-\phi(t)}\leq C.
\]
In view of the two above estimates (\ref{NG}) yields $f(t)\leq C\left( f(0)e^{-Mt}+C\right)$.
\end{proof}
\noindent
\begin{proof} [Proof of Theorem \ref{glo}]
First we argue that for any $p \in [1, \infty)$ it holds
\begin{equation}
|u|_p (t) \le C(u_0, v_0, p, \alpha).
\label{2.0.1}
\end{equation}
In the case of subcritical diffusion we have it for any initial mass by Corollary \ref{cor2.2}. In the critical diffusion
case we need the assumption on smallness of $\chi$ and the initial mass $M$. Then we are in a position to apply
Lemma \ref{lem2.5} and obtain global in time bounds for $|u+1|_{LlogL}$, which in turn gives $|u+1|_{3} (t)\le C$
via Lemma \ref{lem2.8}. This allows us by Corollary \ref{cor2.3} to obtain \eqref{2.0.1}. In order to perform the final step, i.e. reach the bound $|u(t)|_\infty \le C$, we resort to Proposition \ref{glo.5}. We make the following choices in the setting of this theorem, in line with its assumptions
\begin{itemize}
\item[(1)] $D(x,t,s) : = a(s)$, then by definition of $a(s)$, $ D(x,t,s) \in C^\infty ([0,1] \times [0, \infty)^2)$ and it holds: $ D(x,t,s) \ge (1+u)^{-1} \ge (2s)^{-1}$ for $s \ge 1$. Therefore assumption (1) of Proposition \ref{glo.5} holds with $\delta = \frac{1}{2}, \: m=0, \: s_0 = 1$.
\item[(2)] $f := uv_x$. By Proposition \ref{prop1.2} $uv_x \in C^1( (0,T) \times [0,1] )$ and as [KS] boundary condition is Neumann's one, $uv_x(0)=uv_x(1)=0$. Moreover by Proposition \ref{prop1.3} and \eqref{2.0.1} one has $uv_x \in L^\infty(L^{q_1})$ for any $q_1 < \infty$; fix $q_1=17$.
\item[(3)] $u \ge 0$ with the above choices solves [KS] and fixing $p_0 = 2$, $u \in L^\infty(L^{p_0})$ because of \eqref{2.0.1}.
\item[(4)] $17 = q_1 > 3, \: 2=p_0 >1$.
\end{itemize}
Thereby assumptions of Proposition \ref{glo.5} are fulfilled and by its thesis we have $|u+1|_\infty (t) \le C$
\end{proof}
\begin{rem}
Theorem \ref{glo} can be generalized to hold for the system considered in \cite{clKS}, i.e.
\begin{multline}
\begin{cases}
u_t = (a (u) u_x - u v_x)_x \quad in \; (0, \infty) \times (0,1) \\
\varepsilon v_t = D v_{xx}  + u - M +\gamma v \quad in \; (0, T) \times (0,1)  \\
a(u) u_x =  v_x = 0  \quad on \; (0, T) \times \{0;1\}  \\
(u, v) (0) = (u_0, v_0)
\end{cases}
\end{multline}
where $\varepsilon, \gamma, D$ are nonnegative.
\end{rem}
\section{Finite-time blowup}
\noindent
In this section we analyze the opposite situation to that of the previous section, i.e. a possibility of a finite-time
blowup of $|u(t)|_\infty$, $u$ being the solution to [KS]. Our method is a slight extension of results in
\cite{clKS, clSP}. Actually we are going to modify a method in \cite{clKS} in the spirit of \cite[Theorem 10]{clSP}.
The reason for it is that we want to include a wider class of nonlinearities in our result than the ones covered in \cite{clKS}.
For a simplicity of presentation we assume $\chi=1$.
\begin{theo}
Let $\chi=1$ in [KS].For any diffusion $a \in C [0, \infty) \cap L^1 (\er_+) \cap L^\infty (\er_+)$ there exist: such a small $\varepsilon>0$ and $(u_0, v_0)$ with initial mass $M=\int_0^1 u_0(x)dx$ large enough such that a solution $u$ to [KS] emanating from initial data $(u_0,v_0)$ blows up in a finite time.
\label{lem3.4}
\end{theo}
\noindent
We postpone the proof of this result until several technical propositions are proven.
Let us begin with some definitions
\begin{mydef} Let $L, U, V $ be
\begin{itemize}
\item [(i)]$ L (t) : =\frac{1}{q} \int_0^1 |U (t,z)|^q dz$, where $U (t, x) : = \int_0^x u(t, z) dz \;$ for $ x \in [0, 1]$ and $q>2$;
\item[(ii)] $V (t, x) := \int_0^x v(t, z) dz \;$ for $ x \in (0, 1)$.
\end{itemize}
\label{def3.1}
\end{mydef}
\noindent
For a function $B$ (that it will be specified later), $M$- the initial mass of [KS] and $q>2$, let us formally introduce function $A_{B, q} (L) $ as follows
\begin{mydef}
$A_{B, q}(L)  :=  (q-1) B^{\frac{2}{q}} (M) \Big[ \frac{M^{q+1}}{q+1} \Big]^{\frac{q-2}{q}}  \beta^{\frac{q-2}{q}} \Big(\frac{M^{q+1}}{L q(q+1)} \Big)  +  M L [1 +\frac{\varepsilon M^{q-1}}{4q}] -\frac{M^{q+1}}{q(q+1)}
 \quad$
 for $\,\beta (x) : = \frac{B(x)}{x}$.
 \label{def3.2}
\end{mydef}
Next we state the following technical result
\begin{prop}
If there exists a concave real function $B$, for which $\lim_{x \rightarrow \infty} \frac{B(x)}{x} = 0$ and
$0 \le r \int_r^\infty a(s)ds \le B(r)$ where $a(s): [0,\infty)\rightarrow \er_+$ is the diffusion in [KS], then the following differential inequality holds
\begin{equation}
\frac{d}{dt} \left(L + \lambda + \frac{M^2}{2} \right) (t) \le A \left[ \left(L + \lambda + \frac{M^2}{2} \right) \right],
\label{3.3.1}
\end{equation}
where L complies with Definition \ref{def3.1} and $\lambda$ denotes the Liapunov functional associated to [KS] by Proposition \ref{prop1.2}.
\label{prop3.3}
\end{prop}
\begin{proof}
By the nonnegativity of $u$ we have: $L (t) = \frac{1}{q} \int_0^1 |U (t,z)|^q = \frac{1}{q} \int_0^1 U (t,z)^q $ and by differentiation
\begin{equation}
\frac{d}{dt} L (t) = \int_0^1 U (t,z)^{q-1} U_t (t, z) dz =  \int_0^1 U (t,z)^{q-1} [a(u) u_z - uv_z] (t,z) dz.
\label{3.3.2}
\end{equation}
Denoting $A(s) : = - \int_s^\infty a (z) dz$, we see that by \eqref{3.3.2} it holds
\begin{equation}
\frac{d}{dt} L (t) =  \int_0^1 U^{q-1}[A(u)]_z -  \int_0^1 U^{q-1}  u[\varepsilon V_t - U + Mz].
\label{3.3.3}
\end{equation}
Integrating by parts we have
\begin{multline}
\frac{d}{dt} L (t) = - (q-1) \int_0^1 U^{q-2} u A(u) + [U^{q-1} A(u)]\Big|_0^1 +   \int_0^1 U^{q-1}  u[ U - Mz - \varepsilon V_t]  =\\ (q-1) \int_0^1 U^{q-2} [ -u A(u)] + M^{q-1} A(u (1,t))  + \frac{1}{q+1} \int_0^1 (U^{q+1})_z - \frac{1}{q} \int_0^1 (U^{q})_z [Mz + \varepsilon V_t].
\label{3.3.4}
\end{multline}
Observe that $A(u (1,t)) \le 0$ because of the nonnegativity of $a$ as well as the definition of $A$. Next
\[
-u A(u) = u \int_u^\infty a(s) ds \le B(u),
\]
and we arrive at
\begin{multline}
\frac{d}{dt} L (t) \le  (q-1) \int_0^1 U^{q-2} B(u)  + \frac{1}{q+1}M^{q+1} + \frac{1}{q} \int_0^1 U^{q} M - \frac{M}{q} M^q +  \frac{\varepsilon }{q} \left[ \left( \int U^q v_t \right) -  M^q V_t  (1) \right] \\
\le   (q-1) \int_0^1 U^{q-2} B(u) + ML - \frac{M^{q+1}}{q(q+1)} +  \frac{\varepsilon }{q^\frac{1}{2}} M^\frac{q}{2} L^\frac{1}{2} |v_t|_2.
\label{3.3.5}
\end{multline}
Let us focus on the term $ \int_0^1 U^{q-2} B(u)$. We use the estimates from \cite[Theorem 10]{clSP} involving Jensen's inequality for probabilistic measures $\frac{B(u) dx}{\int B(u)}$, $dx$ and  $\frac{U^q dx}{q L}$. We present them for the sake of completeness.
\begin{multline}
 \int_0^1 U^{q-2} B(u)dx =\left( \int_0^1 B(u) dx\right) \int_0^1 (U^{q})^\frac{q-2}{q} \frac{B(u) dx}{\int_0^1 B(u)dx} \\ \le \left( \int_0^1 B(u) dx\right)  \left(\int_0^1 U^{q} \frac{B(u) dx}{\int_0^1 B(u)dx} \right)^\frac{q-2}{q} =
 \left( \int_0^1 B(u) dx\right)^\frac{2}{q} (qL)^\frac{q-2}{q}   \left(\int_0^1 B(u) \frac{U^q dx}{q L} \right)^\frac{q-2}{q} \\ \le
 (qL)^\frac{q-2}{q}  \left[ B \left( \int_0^1 u \right) \right]^\frac{2}{q}  \left[ B \left(  \int_0^1 \frac{u U^q dx}{q L} \right) \right]^\frac{q-2}{q}.
 \label{3.3.6}
\end{multline}
Recalling that $\beta (x): = \frac{B(x)}{x}$ and in virtue of $ B \left(  \int_0^1 \frac{u U^q dx}{q L} \right)  = B \left(\frac{ U^{q+1} |^1_0}{(q+1) q L} \right)$, inequalities  \eqref{3.3.5}, \eqref{3.3.6} yield
\begin{multline}
\frac{d}{dt} L(t) \le  (q-1) B^{\frac{2}{q}} (M) \left[ \frac{M^{q+1}}{q+1} \right]^{\frac{q-2}{q}}  \beta^{\frac{q-2}{q}} \left(\frac{M^{q+1}}{L q(q+1)} \right)  + M L -\frac{M^{q+1}}{q(q+1)} +
\frac{\varepsilon}{q^\frac{1}{2}} M^{\frac{q}{2}} L^{\frac{1}{2}} |v_t|_2.
\label{3.3.7}
\end{multline}
Adding \eqref{3.3.7} and (\ref{Logr}) we arrive at
\[
\frac{d}{dt}  (L + \lambda)(t) \le A_{B, q}(L) - \varepsilon \left[|v_t|_2 + \frac{1}{2q^\frac{1}{2}} M^{\frac{q}{2}} L^\frac{1}{2} \right]^2 \le A_{B, q}(L).
\]
Hence
\begin{equation}
\frac{d}{dt}  \left(L + \lambda + \frac{M^2}{2}\right)(t) \le A_{B, q}(L).
\label{3.3.8}
\end{equation}
In view of the bound $ \lambda (t) \ge - \frac{M^2}{2}$, see (\ref{Logr}), \eqref{3.3.8} implies the claim, provided
the function $A_{B, q}: \er_+ \rightarrow \er$ is nondecreasing. Observe that it suffices to check that the term
$\beta^{\frac{q-2}{q}} \left(\frac{M^{q+1}}{L q(q+1)} \right)$ is nondecreasing. Its monotonicity corresponds to the one of $\beta(x^{-1})$, which is nondecreasing if and only if $\beta(x) = \frac{B(x)}{x}$ is nonincreasing. Assume the contrary, i.e. there are $x_0 > y_0 > 0$, such that $\frac{B(x_0)}{x_0} > \frac{B(y_0)}{y_0}$. Consequently, the graph of $B$ over the interval $[y_0, x_0]$ must remain above the interval in the plane joining points $(x_0, B(x_0))$ and  $(y_0, B(y_0)) $ (denote it by $I$), because of concavity of $B$. Hence tangent to the graph of $B$ at $(y_0, B(y_0)) $ must be at least as steep as the straight line containing $I$, which forces the graph of $B$ to lie below zero for some strictly positive values.
This violates the nonnegativity of $B$.
\end{proof}
\noindent
In order to make use of Proposition \ref{prop3.3} we need the following result, which can be found in \cite[Lemma 3.1]{clGlob_vs_blow}.
\begin{prop}
Let $a \in C^1(0, \infty) \cap L^1 (0, \infty)$, $a \ge 0$ and $\sup_{r \in (0,1)} r \int_r^\infty a(s)ds < \infty$, then there exists a concave function $B\in C[0,\infty)$ such that $B(r) \ge r \int_r^\infty a(s)ds $ and $\lim_{r \rightarrow \infty} \frac{B(r)}{r} = 0$.
\label{prop3.5}
\end{prop}
\begin{rem}
Using the function constructed in \cite[Lemma 3.1]{clGlob_vs_blow} automatically allows us to justify that $\beta(x)$ is nonincreasing. However we have already shown it  for general $B$.
\end{rem}
\begin{proof} [Proof of Theorem \ref{lem3.4}]
By Proposition \ref{prop3.5} assumptions of Proposition \ref{prop3.3} are fulfilled, hence we obtain
\begin{equation}
\frac{d}{dt}  \left(L + \lambda + \frac{M^2}{2}\right)(t) \le A_{B, q}\left[ \left(L + \lambda + \frac{M^2}{2}\right)(t) \right].
\label{3.4.0}
\end{equation}
We aim now at showing that for some $u_0, v_0$ one has
\begin{equation}
 A_{B, q}\left[ \left(L + \lambda + \frac{M^2}{2}\right)(0) \right] < 0.
\label{3.4.1}
\end{equation}
To this end recall that by definition (compare Proposition \ref{prop1.2})
\begin{equation}
\lambda (u,v) := \int \left[b(u) -uv + \frac{1}{2} |v_x|^2\right],
\label{3.4.2}
\end{equation}
where $b: \er_+ \rightarrow \er$ is a function satisfying: $ b'' (s) = \frac{a(s)}{s}, \;  b' (1) = 0, \; b(1) = 0$, therefore one can set
\begin{equation}
b(x) = \left\{ \begin{array}{clcr}
\int_x^1 \left( \int_r^1 \frac{a(s)}{s} ds \right) dr & for \; x \le 1 , \\
\int_1^x \left( \int_1^r \frac{a(s)}{s} ds \right) dr & for \; x \ge 1 ,
\end{array} \right.
\label{3.4.4}
\end{equation}
which implies the following bound on $b$,
\begin{equation}
b(x) \le  |a|_\infty \left| \int_x^1 \left( \int_r^1 \frac{ds}{s} \right) dr \right| = C |1-x + xlnx|  \le \left\{ \begin{array}{clcr}
 C (1+ x)  & for \; x \le 1  \\
 C (1+ x^2) & for \; x \ge 1
\end{array} \right\}
\le  C (1+ x^2).
\label{3.4.45}
\end{equation}
Using the inequality $|v|_\infty \le K |v|_{1,2}$ one can merge \eqref{3.4.2} with \eqref{3.4.45} into
\begin{equation}
\lambda (u,v) + \frac{M^2}{2} \le C  \int u^2 + |v|_\infty \int u + \frac{1}{2} |v|_{1,2}^2  + \frac{M^2}{2} + C \le C  \int u^2 + K |v|_{1,2} \int u + \frac{1}{2} |v|_{1,2}^2  + \frac{M^2}{2} + C.
\label{3.4.6}
\end{equation}
Since it is satisfactory to give any examples of $u_0, v_0$ leading to the validity of \eqref{3.4.1}, for $M>1$ we pick the following initial data
\begin{equation}
 \begin{array}{crcccl}
u_0 & = & 2M^3\left( x + \frac{1}{M} - 1 \right) \id_{ [1 - \frac{1}{M}, 1] } (x),  \\
v_0 & = &  Mx -  \frac{M}{2},
\end{array}
\label{3.4.7}
\end{equation}
which yields
\begin{equation}
\int u_0 = M, \quad \int u_0^2 = \frac{4}{3} M^3, \quad |v_0|_{1,2}^2 = \frac{13}{12} M^2,
\label{3.4.8a}
\end{equation}
\begin{equation}
L(0) := \frac{1}{q} \int_0^1  \left[\int_0^x u_0 (z) dz\right]^q dx = \frac{2^q M^{q-1} }{q(2q+1)} .
\label{3.4.8b}
\end{equation}
In view of \eqref{3.4.8b}, \eqref{3.4.8a}, \eqref{3.4.6} one estimates
\begin{equation}
 \frac{2^q M^{q-1} }{2q(q+1)}  \le L(0) + \lambda (u_0,v_0) + \frac{M^2}{2} \le  \frac{2^q M^{q-1} }{(q+1)q} + CM^3 + C.
\label{3.4.85}
\end{equation}
Recall that for $B$ it holds $\lim_{x \rightarrow \infty} \frac{B(x)}{x} = 0$. Let us choose $M_0$ such that for any $M \ge M_0$
\begin{equation}
B(M) \le M, \quad B \left(\frac{M^{q+1}}{\left( L(0)  +\lambda (u_0,v_0) + \frac{M^2}{2} \right)  q(q+1)} \right) \le c \frac{M^{q+1}}{\left( L(0) + \lambda (u_0,v_0) + \frac{M^2}{2} \right)  q(q+1)},
\label{3.4.87}
\end{equation}
while the latter is possible for fixed $q>2$ in view of \eqref{3.4.85}.  Let us now calculate, using the definition of $A$ as well as \eqref{3.4.87}
\[
A_{B, q}\left[ \left(L + \lambda + \frac{M^2}{2}\right)(0) \right] =
\]
\[
(q-1) q^{q-2}B^{\frac{2}{q}} (M) \left[ L(0)+  \lambda (u_0,v_0) + \frac{M^2}{2} \right]^{\frac{q-2}{q}}  B^{\frac{q-2}{q}} \left(\frac{M^{q+1}}{\left( L(0)  + \lambda (u_0,v_0) + \frac{M^2}{2} \right) q(q+1)} \right)
\]
\[
+  M\left[1 +\frac{\varepsilon M^{q-1}}{4q}\right] \left[L(0)+ \lambda (u_0,v_0) + \frac{M^2}{2}\right]  -\frac{M^{q+1}}{q(q+1)} \le
\]
\[
c(q) M^{\frac{2}{q}}  M^{\frac{(q+1)(q-2)}{q}} + cM \left[ \frac{2^q M^{q-1} }{2q(q+1)} + CM^3 + C \right] -\frac{M^{q+1}}{q(q+1)} \le
\]
\[
C(q) [M^q+ M^4 + 1] - \frac{M^{q+1}}{q(q+1)},
\]
where we chose $\varepsilon:=\frac{1}{M^{q-1}}$, which gives us the smallness condition on $\varepsilon$ already mentioned
in the Theorem \ref{lem3.4}. Therefore assuming $q>4$, the bound  \eqref{3.4.1} holds for $M$ large enough. Knowing this, we
conclude in a following standard way that $|u|_\infty$ blows up in a finite time. First assume that $T_{\max} = \infty$,
then by \eqref{3.4.1} and \eqref{3.4.0} we obtain that for finite time $t>0$, $L(t) $ becomes negative, which is absurd. By Proposition \ref{prop1.1} this implies $T_{\max} < \infty$ and a finite-time blowup.
\end{proof}

\section{Conclusions}
\noindent
The present paper has been aimed as a first step in the studies of correspondence between qualitative behaviour of
solutions to one-dimensional quasilinear parabolic-elliptic and fully parabolic Keller-Segel systems. As we mentioned
in the introduction, this question is worth of studies, since in higher dimensions and biologically more relevant cases,
rigorous results are (almost always) presented in a parabolic-elliptic case, while this is a fully parabolic one which
forms an original model. We focused on the question whether in a one-dimensional fully parabolic case, like in the 
corresponding parabolic-elliptic one,
there is no critical nonlinearity. It is known that such a phenomenon takes place in the latter (see \cite{clSP}, \cite{clGlob_vs_blow}).
We did not succeed in  answering this question, however we managed to make a first step in studying this problem. On the
one hand we fully proved global-in-time existence of solutions in subcritical cases, on the other we proved finite-time
blowup of solutions when nonlinear diffusion is integrable at infinity.

Moreover, our studies of the case when a nonlinear diffusion is of the form $\frac{1}{u+1}$ led us to considerations of  
nonlinear functional inequalities of the type (\ref{sehrwichtig}). One could expect that proving this inequality is a proper way
of achieving the global-in-time existence of solutions to [KS] without restriction on the size of initial mass.
As stated in Remark \ref{uwaga}, it is false. Here we give a counterexample.
 \begin{lem}
Let a class of functions ${\cal A}$ be defined in a following way. We say that $m\in W^{1,2}(0,1)\in {\cal A}$ if 
\begin{equation}\label{wiez}
\int_0^1 e^m dx=M.
\end{equation}
Then for any small enough $\delta>0$, independently on the choice of a continuous function $h_\delta: [0,\infty)\rightarrow \er$,
the inequality 
\begin{equation}\label{sehrwichtig2}
\int e^{2m} < \delta M^2\int |m_x|^2 + h_\delta \left( \int e^m \right)
\end{equation}
does not hold in all the class ${\cal A}$. 
 \label{nieprawda}
 \end{lem}
 \begin{proof}
Fix an arbitrary $\delta > 0$ and a continuous function $h_\delta$. We define a set of functions 
$m_\varepsilon (x)\in W^{1,2}(0,1)$ satisfying (\ref{wiez}) in a following way
\[
e^{m_\varepsilon (x)} = \frac{\varepsilon(1+\varepsilon)M}{(x+\varepsilon)^2}.
\]
Then, on the one hand (\ref{wiez}) is satisfied for each $\varepsilon>0$ and on the other hand
\[
\delta M^2\int |m_x|^2 + h_\delta \left( \int e^m \right) = \frac{4 \delta M^2 }{\varepsilon (1+\varepsilon)} + h_\delta \left(M\right) \le \frac{4 \delta M^2 }{\varepsilon} + C,
\]
where $C$ is independent on $\varepsilon$. Moreover,
\[
\int e^{2m_\varepsilon} = \frac{M^2}{3} \frac{(1+ \varepsilon)^3 - \varepsilon^3 }{\varepsilon (1+\varepsilon)}. 
\]
Next notice that taking $\delta\leq 1/24$ the right-hand side of (\ref{sehrwichtig2}) is less than 
$\frac{M^2 }{6\varepsilon} + C$ while the left-hand side of (\ref{sehrwichtig2}) behaves like 
$\frac{M^2}{3\varepsilon}$  for $\varepsilon$ close to $0$.
\end{proof}

Finally let us mention that our finite-time blowup result holds true only for large initial masses and small factors
$\varepsilon>0$ by the time derivative of $v$. It is an interesting open problem if one can, like in the parabolic-elliptic case, achieve this result for any size of the initial mass.
Another challenging problem is: if the finite-time explosion of the solution can happen for $\varepsilon>0$ large.
\vskip 5mm
\noindent
{\bf Acknowledgement.} The first author is thankful to prof. K. Oleszkiewicz for his valuable remark.




\noindent

\end{document}